\documentclass[12pt]{article}
\usepackage{amsmath,amssymb,amsfonts}
\newcommand{\eps}{\epsilon}
\newcommand{\Z}{\mathbb Z}
\renewcommand{\H}{\mathbb H}
\newcommand{\E}{\mathbb E}
\newcommand{\C}{\mathbb C}
\newcommand{\R}{\mathbb R}
\newcommand{\T}{\mathbb T}
\renewcommand{\Re}{\text{Re}}
\renewcommand{\Im}{\text{Im}}
\newtheorem{thm}{Theorem}
\newtheorem{prop}{Proposition}
\newtheorem{lemma}{Lemma}

\newenvironment{proof}{
    \noindent{\bf Proof:} \hspace*{1em}}{
    \hspace*{\fill} $\square$\medskip }
\newcommand{\PSbox}[3]{\mbox{\rule{0in}{#3}\hspace{#2}\includegraphics{#1}}}

\begin{document}
\title{An introduction to the dimer model}
\author{Richard Kenyon
\thanks{Laboratoire de Math\'ematiques,
UMR 8628 du CNRS, B\^at 425, Universit\'e Paris-Sud,
91405 Orsay, France.}}
\maketitle

\section{Introduction}
A {\bf perfect matching} of a graph is a subset of edges which covers
every vertex exactly once, that is, for every vertex there is exactly
one edge in the set with that vertex as endpoint. 
The {\bf dimer model} is the study of the set of perfect matchings of a
(possibly infinite) graph. The most well-known example is
when the graph is $\Z^2$, for which perfect matchings are equivalent
(via a simple duality) to {\bf domino tilings}, that is, tilings of
the plane with $2\times 1$ and $1\times 2$ rectangles.

In the first three sections we study domino tilings of the plane
and of finite polygonal regions, or equivalently, perfect matchings
on $\Z^2$ and subgraphs of $\Z^2$.

In the last two sections we study the FK-percolation model and the dimer
model on a more general family of planar graphs.

\section{The number of domino tilings of a chessboard}
\begin{figure}[htbp]
\vskip3in
\PSbox{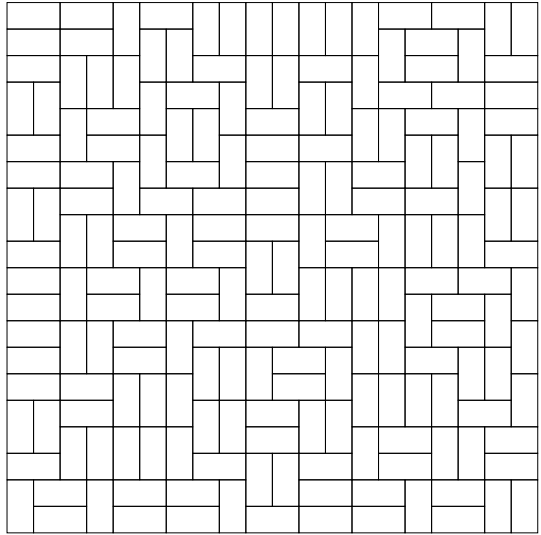}{0in}{0in}
\caption{Random tiling of a square}
\end{figure}

A famous result of Kasteleyn \cite{Kast1} and Temperley and Fisher
\cite{FT} counts the number
of domino tilings of a chessboard (or any other rectangular region).
In this section we explain Kasteleyn's proof.

\subsection{Combinatorics}
Let $R$ be the region bounded by
a simple closed polygonal curve in $\Z^2$. 
A domino tiling of $R$ corresponds to a perfect matching of $G$,
the dual graph of $R$: $G$ has a vertex for each lattice square in $R$,
with two vertices adjacent if and only if the corresponding lattice
squares share an edge.

\begin{thm}\label{1}[Kasteleyn, 1961] 
The number of domino 
tilings of $R$ is $\sqrt{|\det K|}$, where $K$ is the weighted adjacency
matrix of the graph $G$, with horizontal edges weighted
$1$ and vertical edges weighted $i=\sqrt{-1}$.
\end{thm}

\begin{figure}[htbp]
\vskip2in
\PSbox{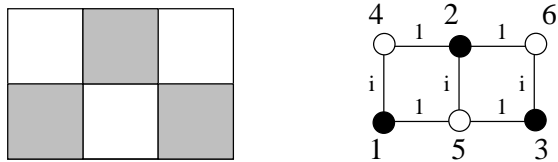}{0in}{0in}
\caption{$2\times 3$ rectangle and dual graph.\label{2by3}}
\end{figure}

For example for the $2\times 3$ region in Figure \ref{2by3},
the matrix $K$ is
$$K=\left(\begin{array}{cccccc}0&0&0&i&1&0\\0&0&0&1&i&1\\
0&0&0&0&1&i\\i&1&0&0&0&0\\1&i&1&0&0&0\\0&1&i&0&0&0
\end{array}\right),$$
whose determinant has absolute value $9=3^2$.

\begin{proof}
Since the graph $G$ is bipartite (one 
can color vertices black and white so that black vertices are
only adjacent to white vertices and vice versa)
$K$ can be written $K=\left(\begin{array}{cc}0&A\\A^t&0\end{array}\right)$,
so we must evaluate 
$$\det A = \sum_{\sigma\in S_n} \text{sgn}\sigma 
k_{1\sigma(1)}\dots k_{n\sigma(n)}$$
Each term in this sum is zero unless black vertex $\sigma(i)$
is adjacent to white vertex $i$ for each $i$.
Therefore we have exactly one term for each matching of $G$.
It suffices to show that two nonzero terms have the same sign.

Suppose we have two matchings; draw them one on top of the other.
Their superposition is a union of (doubled) edges and disjoint 
cycles. 
One can transform the first matching into the second by 
``rotating'' each cycle in turn.
In particular it suffices to show that if two matchings
differ only along a single cycle, the corresponding terms
in the determinant have the same sign.

This is supplied by the following lemma,
which is easily proved by induction.
\end{proof}

\begin{lemma}\label{pathcomb}
Let $C=\{v_0,\dots,v_{2k-1},v_{2k}=v_0\}$ be a simple closed curve in $\Z^2$.
Let $m_1$ be the product of the weights on edges
$v_0v_1, v_2v_3,\dots$ and $m_2$ the product of the weights on the
remaining edges $v_1v_2,\dots$ (recall that vertical edges are weighted $i$
and horizontal edges $1$). Then $m_1=(-1)^{n+k+1} m_2$, where $n$ is the 
number of vertices strictly enclosed by $c$ and $k$
is $1/2$ of the length of $c$.
\end{lemma}

\subsection{Rectangles}
Suppose the graph $G$ is an $m$ by $n$ rectangle, with vertices
$\{1,2,\dots,m\}\times\{1,2,\dots,n\}$.
To compute the determinant of its adjacency matrix,
we compute its eigenvectors and their eigenvalues.

For $j,k\in\Z$ 
define $$f(x,y)=\sin\frac{\pi jx}{m+1}\sin\frac{\pi ky}{n+1}.$$
It is easy to check that $f$ is an eigenvector of $K$,
with eigenvalue 
$$2\cos\frac{\pi j}{m+1}+2i\cos\frac{\pi k}{n+1}.$$
We didn't pull this formula out of a hat: since $G$ is a
``product'' of two line graphs, the eigenvectors are the product of the
eigenvectors of the line graphs individually.

As $j,k$ vary over integers in $[1,m]\times[1,n]$,
the $f$ form an orthogonal basis of eigenvectors.
Therefore the determinant of $K$ is
$$\det K = \prod_{j=1}^m\prod_{k=1}^n 
2\cos\frac{\pi j}{m+1}+2i\cos\frac{\pi k}{n+1}.$$

Evaluating $\sqrt{|\det K|}$ for $m=n=8$ gives $12988816$ tilings of a
chessboard.

One may show that for $m,n$ both large the number of tilings is
$\exp(Gmn/\pi+O(m+n)),$ where $G=1-\frac1{3^2}+\frac1{5^2}-\frac1{7^2}+\dots$
is Catalan's constant.

\subsection{Tori}
It will be useful shortly to compute the number of tilings for a graph
on a torus as well (e.g. a rectangle with its opposite sides identified).
Kasteleyn showed that this could be accomplished 
by a linear combination of four determinants.
In essence these four determinants correspond to the four
inequivalent ``discrete spin structures" on the graph $G$
embedded on the torus. 
Rather than go into details, we just note that the above proof
fails on a torus since two matchings may differ on a loop
which goes around the torus (is not null-homologous).
The change in sign in the determinant 
from one matching to the next is then a function
of the homology class in $H_1(\T,\Z/2\Z)$ of the loop.
The result is that (for $m,n$ both even)
the number of tilings of an $m\times n$ torus is
$$Z_{m,n}=\frac12(-P_{00}+P_{01}+P_{10}+P_{11}),$$
where 
$$P_{\sigma\tau}=\prod_{z^m=(-1)^\sigma}\prod_{w^n=(-1)^\tau}
\left(z+\frac1z+iw+\frac{i}{w}\right).$$
Note that actually $P_{00}=0$ in this case.
For details see \cite{Kast1} or \cite{CKP}. 

One can show that 
$$\lim_{m,n\to\infty}\frac1{mn}\log Z_{m,n}=\frac{G}{\pi}.$$

\subsection{Inverse Kasteleyn matrix}
We are primarily
interested in the {\bf uniform measure} on perfect matchings of a
graph $G$.

A very useful consequence of Kasteleyn's counting theorem is
\begin{thm}[\cite{localstats}]\label{ls}
Let $T=\{(w_1,b_1),\dots,(w_k,b_k)\}$ be a subset of dimers,
(with the $j$th dimer covering white vertex $w_j$ and black vertex $b_j$).
The probability that all dimers in $T$ appear in a
uniformly chosen matching is
$$\Pr(T) = |\det(K^{-1}(w_i,b_j))_{1\leq i,j\leq k}|.$$
\end{thm}
The advantage of this result is that the determinant for the joint
probability of $k$ dimers is of size $k$
{\it independently of the size of the matrix $K$}.
In particular to compute the probability that a single
edge appears, one just needs a single element of $K^{-1}$.

If we can compute asymptotics of $K^{-1}$ for large rectangles
we can then understand easily
the limiting measures on tilings of the plane.
We'll do this later.

A more general class of measures on perfect matchings 
uses graphs with weighted edges.
If edges are weighted $\nu(e)>0$, the weight of a perfect matching
is the product of its edge weights, and the probability measure
we are interested in is that which gives a matching a probability
proportional to its weight (for finite graphs). 
There is a simple generalization of Theorems \ref{1} and \ref{ls} for weighted
graphs.

We should also point out that versions of 
Theorems \ref{1} and \ref{ls} hold for arbitrary planar graphs,
see \cite{Kast2} and \cite{localstats}.

\section{The arctic circle phenomenon}
The order-$n$ Aztec diamond is the region shown.
It was defined and studied in \cite{EKLP, CEP}.
A random tiling is also shown. Why does it not look 
homogeneous, as does the tiling of a square?
\begin{figure}[htbp]
\vskip4in
\PSbox{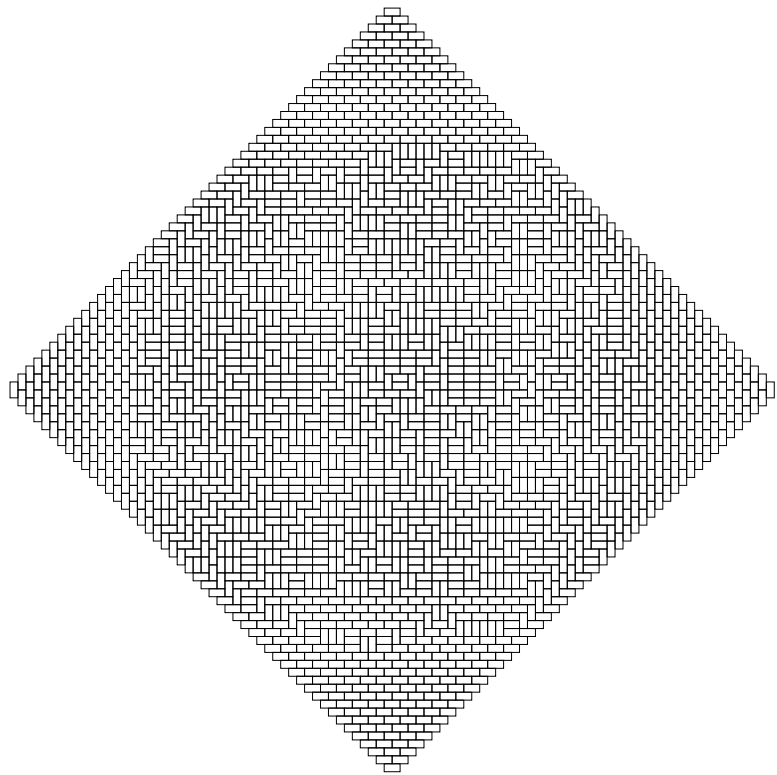}{0in}{0in}
\caption{Random tiling of Aztec diamond}
\end{figure}

Near the corners, you only see tiles of one "type".
The way to understand this phenomenon is via the height function.

\subsection{heights}
A domino tiling can be thought of as an interface in $2+1$ dimensions.
Specifically, to a tiling we associate an
integer-valued function $h$ on the vertices of the region tiled,
as follows. 
Rather than define the height, we define the height 
difference between two adjacent lattice points.
If edge $vw$ is not crossed by a domino, the height difference
$h(w)-h(v)$ is $1$ if the face to the left of $vw$ is black,
and $-1$ otherwise. If the edge is crossed by a domino,
$h(w)-h(v)$ is $-3$ or $+3$ according to whether the face 
to the left is black or white. This defines a ``one-form", that is a
function on the oriented edges of the graph $G$ satisfying
$f(-e)=-f(e)$. 
Moreover the one-form is closed: the sum around any oriented cycle is zero
(since the sum around any face is zero: a domino crosses exactly
one edge of each face). Therefore there is a function
$h$, well defined up to an additive constant, with these edge-differences.

\begin{figure}[htbp]
\vskip3in
\PSbox{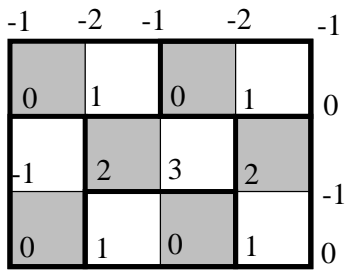}{0in}{0in}
\caption{Heights in a domino tiling.\label{hts}}
\end{figure}

The height function of a tiling gives a map from $\Z^2$ to $\Z$.
One of the principal motivations for studying dominos is 
to study the model of random maps of this sort.

Note that a tiling determines a height function, and vice versa. Also,
the height function on the boundary of the region to be tiled 
is independent of the tiling. 

For the $2n\times 2n$
square, the boundary function is asymptotically flat,
in fact assuming the value at a corner is $0$, it alternates
on edge edge between $0$ and $1$ or $0$ and $-1$: see Figure \ref{hts}.
However for the Aztec diamond of order $n$, it
is linear on each edge:

\begin{figure}[htbp]
\vskip4in
\PSbox{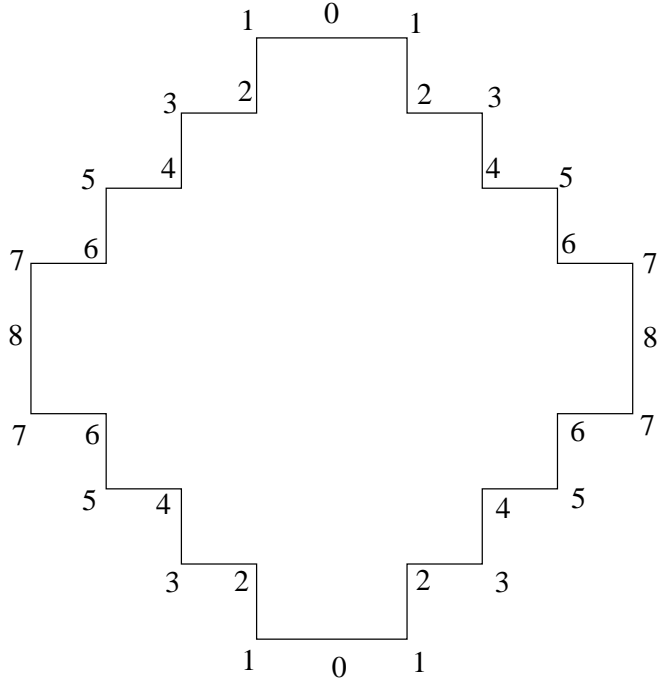}{0in}{0in}
\caption{Heights on boundary of an Aztec diamond.\label{aztechts}}
\end{figure}

For a non-planar region such as a torus, the height might not be
well-defined: it can be well-defined locally, but on a path winding
around the torus it may have a non-trivial period. 
Indeed, for a graph on the $m\times n$ torus the 
height will have a well-defined horizontal period $h_x$ and
vertical period $h_y$, both integers, such that the height
change on any closed path winding once horizontally around the torus
is $h_x$, and similary $h_y$ for vertically-winding paths.

The height of a random tiling on a planar region is a random function 
with those given boundary values. It is easy to 
guess that the height at an interior point of 
a random tiling of a square of side $n$ is 
$o(n)$. However for the Aztec diamond it is not at all clear
what the height of a random tiling will be. In fact it will
typically lie near (within $o(n)$) of its mean value.
Moreover its mean value can be determined by solving a
certain variational principle (minimizing a certain energy).
That is the content of this section.

\subsection{CKP theorem}
The theorem of Cohn-Kenyon-Propp states this in
general: given a Jordan domain $U$
with function $u:\partial U\to\R$, and a sequence of 
domains in $\eps\Z^2$ converging to $U$, with
boundary height functions $h_\eps$ converging to $u$,
that is $\eps h_\eps\to u$, then with probability tending to
$1$ the height of a random
tiling will lie near (within $o(n)$) its mean value,
and the mean value $h_0$ is the unique function $f$
maximizing a certain ``entropy''
$$\mbox{Ent}(f) = \int_U\text{ent}(\partial_x f,\partial_y f)dxdy.$$
Here 
$$\text{ent}(s,t)=
\frac1\pi\left(L(\pi p_a)+L(\pi p_b)+L(\pi p_c)+L(\pi p_d)\right),$$
where 
$L(x)= -\int_0^x\log2\sin t dt$ is Lobachevsky's function,
and $p_a,p_b,p_c,p_d$ are certain probabilities which are determined
by the following equations.
\begin{eqnarray*}2(p_a-p_b)&=&s\\2(p_d-p_c)&=&t\\p_a+p_b+p_c+p_d&=&1\\
\sin\pi p_a\sin\pi p_b&=&\sin\pi p_c\sin\pi p_d.
\end{eqnarray*}

The proof of the first part of this theorem is somewhat standard.
The harder part is the computation of the function $\text{ent}$.

The idea of the proof of the first part is the following.
Each height function is a Lipschitz function with Lipschitz constant $3$,
since $|f(x_1,y_1)-f(x_1+1,y_1)|\leq 3$ and similarly for the $y$-direction.
In fact we have a stronger condition that when
$x_1-x_2,y_1-y_2$ are even, 
\begin{equation}\label{Lip}
|f(x_1,y_1)-f(x_2,y_2)|\leq
2\max\{|x_1-x_2|,|y_1-y_2|\}
\end{equation}
Let $L$ be the class of functions satisfying (\ref{Lip}).
When we scale the $x-$ and
$y$-coordinates by $\eps$, then $\eps f$ is still in $L$.
Note that $L$ is compact (in the uniform topology). 

On this space, the functional
$\text{Ent}$ is upper semicontinuous: its value can only increase at a limit
point.  
This follows from the strict concavity of the function $\text{ent}(s,t)$. 
Concavity also proves unicity of the maximizing function.

Since $L$ is compact, cover it with a finite number of metric balls
$B_f(\delta)$ of radius $\delta$.
We will show that for any $f\in L$, the number of tilings for which
$\eps h$ lies within $\delta$ of $f$ is 
$\exp(\frac1{\eps^2}\text{Ent}(f)(1+o(1)))$.
This determines the size of each $B_f(\delta)$.
Therefore when $\eps$ is small, the number of tilings whose normalized
height function lies close to $f_{\max}$
overwhelmingly dominates the number of other tilings, even combined.

Thus almost all tilings lie near $f_{\max}$.

To count the number of tilings whose height function lies near 
a given function $f\in L$,
triangulate $U$ with a fine mesh, large on the scale of $\eps$ but
with mesh size tending to zero. Since a Lipschitz function is
differentiable almost everywhere, on almost all triangles the
function $f$ is almost linear. A tiling which lies close to $f$
will lie close to a linear function on most triangles.
So it suffices to count the number of tilings whose height function
is close to a linear function. 

However a short technical lemma (essentially subadditivity)
proves that
if $R$ is a triangular region (or other nice region), 
the number of tilings with height function
lying close to a fixed linear function $\ell$ on that region
is $\exp(A\cdot\text{ent}(\ell)\cdot(1+o(1)))$, where $A$ is the area of $R$ and
$\text{ent}(\ell)$ is a function depending only on the slope of $\ell$.

\medskip

The explicit computation of $\text{ent}(s,t)$ comes from 
computing the number of tilings on a torus with a given 
slope $(s,t)$, that is, a given horizontal height change
$[sm]$ going horizontally around the torus, and $[tn]$ vertically.

\begin{figure}[htbp]
\vskip3in
\PSbox{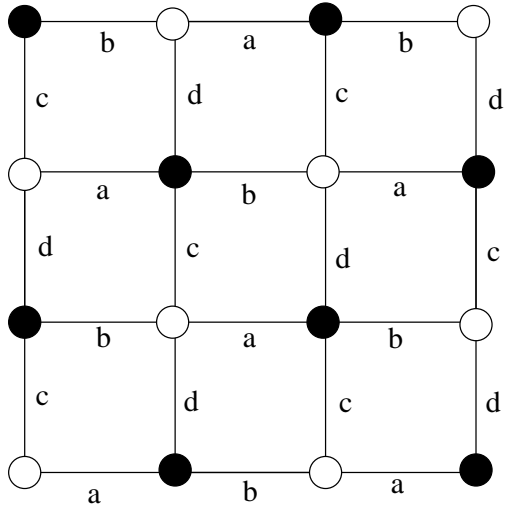}{0in}{0in}
\caption{Weights on torus graph\label{stag}}
\end{figure}

To compute this we use staggered weights $a,b,c,d$ on the edges of $\Z^2$,
as in Figure \ref{stag}.

If we put these weights in the Kasteleyn matrix then its determinant
is the sum over tilings of the weight of a tiling (product
of the edge weights). That is, $\sqrt{|\det K|}$ is the 
{\bf partition function} $Z(a,b,c,d)$ of tilings with these weights.
For a tiling of the torus with these weights, the number of dominos
of type $a,b,c,d$ determines the height change around the torus:
the horizontal height change is $(N_d-N_c)/n$
and the vertical height change is $(N_a-N_b)/m$.
Therefore once we know 
$$Z(a,b,c,d)=\sum_{\text{matchings}} a^{N_a}b^{N_b}c^{N_c}d^{N_d}$$
we have
$$Z(a,1/a,c,1/c)=\sum_{\text{matchings}}a^{mh_y}c^{-nh_x}$$
and we need only extract the appropriate coefficient to get the number
of tilings of given slope.

In \cite{CKP} the asymtotic formula is given
$$\log Z(a,b,c,d)=\frac1{(2\pi i)^2}\int_{S^1}\int_{S^1}
\log(a+b z+cw+dzw)\frac{dz}z
\frac{dw}w,$$
from which the above formula for $\mbox{ent}$ follows.

\subsection{flats}
For the Aztec diamond, the function $f_{\max}$ is not analytic, only piecewise
analytic.
In the four regions outside of the inscribed circle $f_{\max}$ is linear
(as you can approximately see in the example tiling). Inside it is analytic.
On the boundary of the inscribed circle (the {\bf arctic circle}, 
since it separates
the ``frozen region'' from the ``temperate'' zone) the function is
$C^1$ but not $C^2$. In fact it is $C^{1.5}$ except at the four
boundary-edge midpoints.

These frozen regions result from the degeneration of the ellipticity
of the PDE defining $f_{\max}$ at the boundary of the domain of definition.

It is very interesting to study what happens exactly at the boundary.
We won't discuss this here though.

\section{Conformal invariance}
The scaling limit (limit when the lattice spacing tends to zero)
of the height function on domino tilings tends to a nice 
conformally invariant continuous process, the ``massless 
Gaussian free field'', a sort of two-dimensional
version of Brownian motion. 
In this section we discuss the ideas behind the proof.
The original references are \cite{K.confinv,K.gauss}.

\subsection{$K$ as Dirac operator}
\begin{figure}[htbp]
\vskip3in
\PSbox{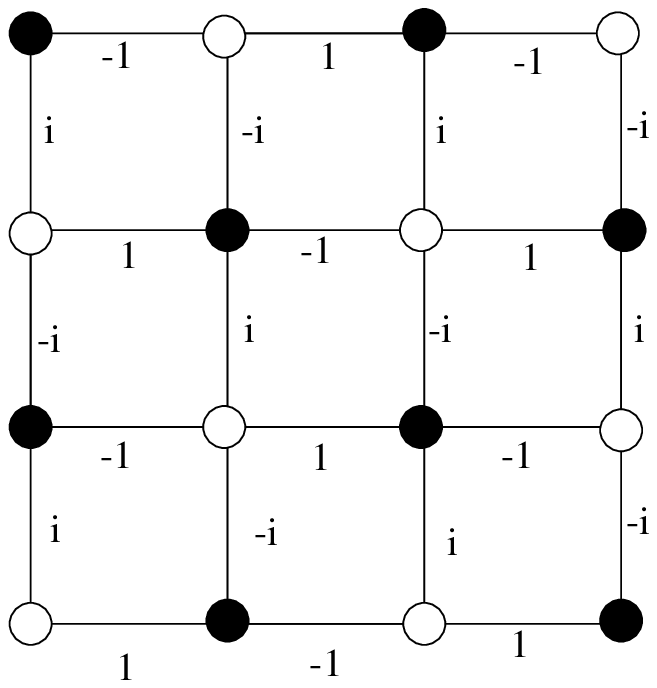}{0in}{0in}
\caption{\label{iwts}}
\end{figure}

We had a lot of choice in definition of the Kasteleyn matrix.
If you go back to the proof, you will see that we could have chosen any
edge weights with the property that the weights on a lattice square
$\alpha,\beta,\gamma,\delta$ satisfy the condition
that $\alpha\gamma/\beta\delta$ is real and negative. It is convenient to 
put weights as in Figure \ref{iwts}.

Lemma \ref{pathcomb} holds for these weights as well, as do
Theorems \ref{1} and \ref{ls}.

Let us now study $K^{-1}$ for these weights.
Because of bipartiteness, $K^{-1}(v_1,v_2)=0$ unless one of
$v_1,v_2$ is black and the other is white. 
Let $B$ denote the black vertices and $W$ the white vertices.
Let $w\in W$. If $f\in\C^B$, then $Kf\equiv 0$ means
$$f(w+1)-f(w-1)+i(f(w+i)-f(w-i))=0$$
that is
``$\partial_x+i\partial_y$''$f=0.$
We say $f\in\C^B$ is {\bf discrete analytic} if $Kf=0$.

Let $B_0,B_1$ be the black vertices whose coordinates
are both even (respectively, both odd).
Similarly let $W_0,W_1$ be the white vertices whose
coordinates are $(1,0)\bmod 2$, respectively $(0,1)\bmod 2$.

\begin{lemma} If $f\in\C^B$ is discrete analytic, then
$f$ is harmonic on $B_0$ and $B_1$ separately,
that is, 
$$4f(b)=f(b+2)+f(b-2)+f(b+2i)+f(b-2i).$$
Moreover, these two harmonic functions are conjugate.
\end{lemma}

\begin{figure}[htbp]
\vskip3in
\PSbox{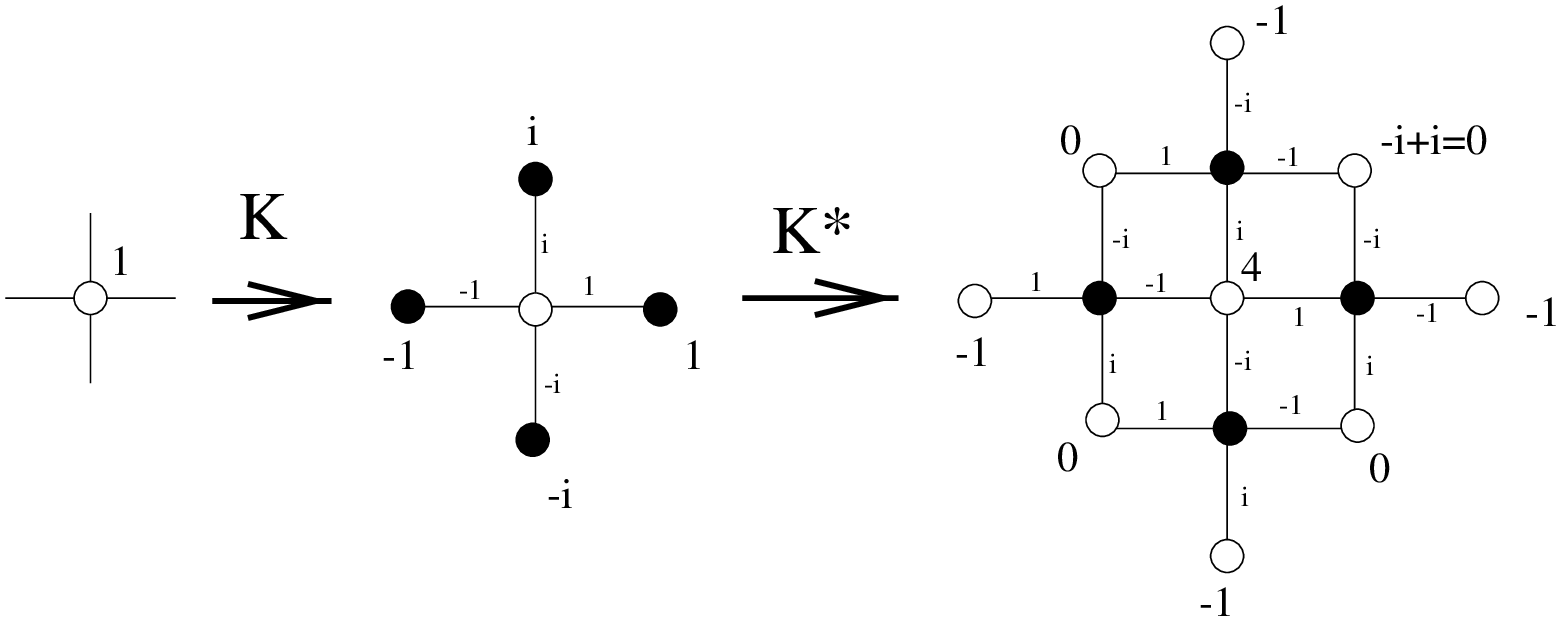}{0in}{0in}
\caption{\label{KK}}
\end{figure}

The proof follows from the observation that 
$K^*K$ is the Laplacian when restricted to any of the sublattices
$B_0,B_1,W_0,W_1$, see Figure \ref{KK}.

As a consequence $K^{-1}$ can be related to the Green's function.

\subsection{$K^{-1}$ and the Green's function}
For a fixed $w$, $K^{-1}(w,b)$ (considered as a function of $b$)
is discrete analytic, except at $w$. This follows
from $KK^{-1}=I$.

In fact $K^{-1}$ can be written in terms of the Green's
function $G$. Recall that $G(u,v)$ is defined to be the
function satisfying $\Delta G(u,v)=\delta_u(v)$ (taking
the Laplacian wrt the second variable). 
For a fixed $w$, 
$$\Delta K^{-1}(w,\cdot) = K^*KK^{-1}(w,\cdot)=
(K^*\delta_w)(\cdot)=\delta_{w+1}-\delta_{w-1}-i(
\delta_{w+i}-\delta_{w-i}).$$
Thus $K^{-1}(w,\cdot)$ is a sum of four Green's functions. 
When restricted to $B_0$, it is the difference of two Greens functions,
$K^{-1}(w,b)=G(w+1,b)-G(w-1,b).$

\subsection{On a rectangle}
We can easily compute $K^{-1}$ for a rectangle, since we already
diagonalized the matrix $K$.
Because $K$ is symmetric, we have the following:
$$K^{-1}(v_1,v_2) = \sum_j f_j(v_1)f_j(v_2)/\lambda_j$$
where $f_j$ are the orthonormalized eigenvectors with eigenvalue $\lambda_j$.

In the limit when the rectangle gets large, and $v_1,v_2$ are far from the
boundaries, we can find
\begin{equation}\label{close}
K^{-1}(v_1,v_2)= \frac1{4\pi^2}\int_0^{2\pi}\int_0^{2\pi}
\frac{e^{i(x\theta+y\phi)}}{2i\sin\theta-2\sin\phi}d\theta d\phi,
\end{equation}
where $(x,y)=v_2-v_1\in\Z^2$.

\begin{figure}[htbp]
\vskip4in
\PSbox{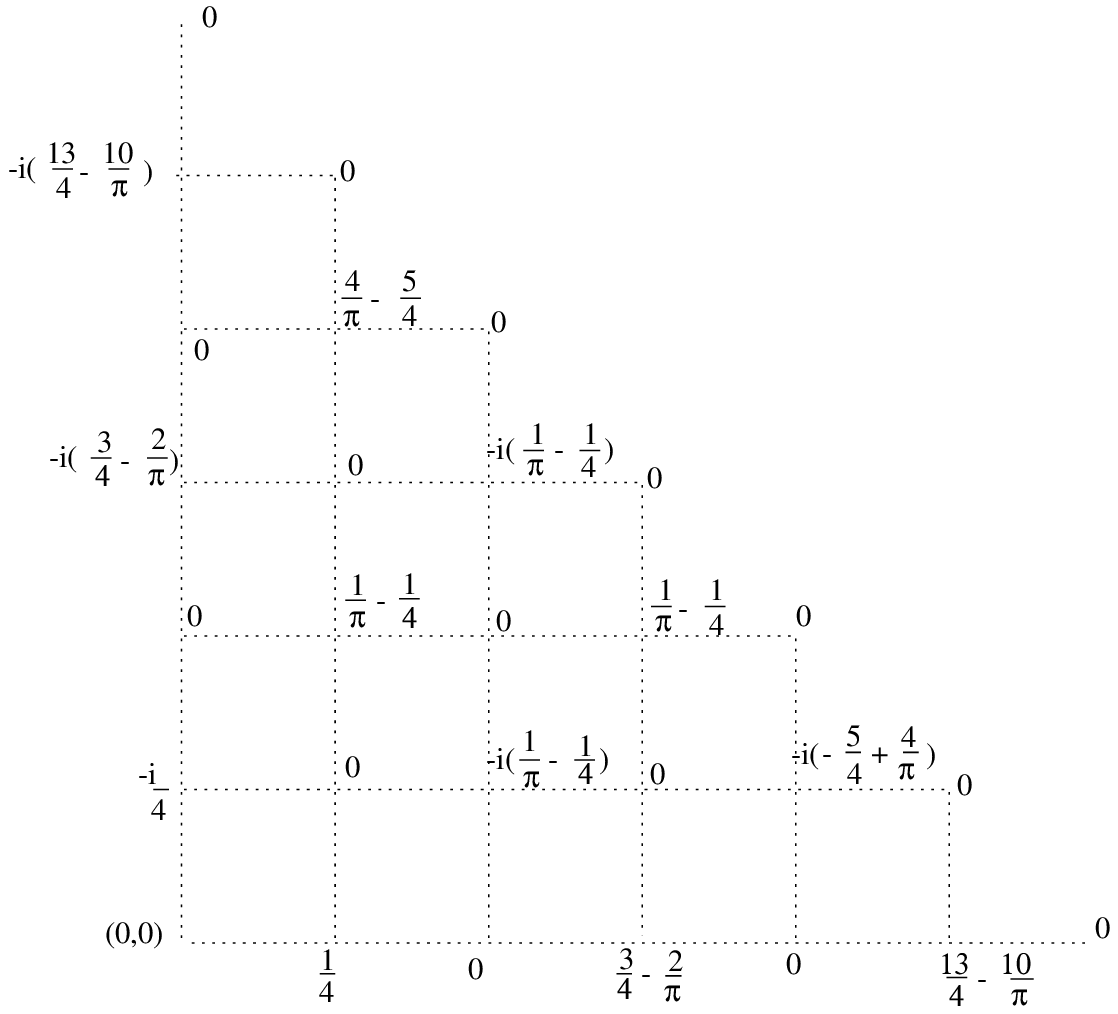}{0in}{0in}
\caption{\label{Kinv}}
\end{figure}

This formula defines, along with Theorem \ref{ls}, a measure on tilings of 
the plane. This measure is the unique weak limit of measures on rectangles.
It is also the unique entropy-maximizing 
translation-invariant measure on tilings \cite{BP}.
Some values of $K^{-1}$ are shown in Figure \ref{Kinv}.

Asymptotically, if $w\in W_0$, we have 
$$K^{-1}(w,b)=\left\{\begin{array}{ll}
\Re\frac1{\pi (b-w)}+O(\frac1{|b-w|^2})&b\in B_0\\
i\Im\frac1{\pi (b-w)}+O(\frac1{|b-w|^2})&b\in B_1.
\end{array}\right.$$
These are reversed if $w\in W_1$.
\subsection{Bounded domains}
We can similarly compute the asymptotics of $K^{-1}$ on bounded domains.

Let $U$ be a Jordan domain. Let $U_{2\eps}$ be a domain
in $2\eps\Z^2$ which approximates $U$. It is a well-known fact that
the discrete Dirichlet Green's function on $U_{2\eps}$
(defined by: $G(x,y)$ is the expected number of passages at $y$ of
a simple random walk started at $x$ and stopped at the boundary,
alternatively, $G(x,y)$ is the inverse of the discrete Laplacian with
Dirichlet boundary conditions)
approximates the continuous Dirichlet Green's function on $U$.

We can use this fact to construct a graph $U_\eps\in\eps\Z^2$
for which the $K^{-1}$ operator converges nicely to the
corresponding continuous object on $U$.

The graph $U_\eps$ has a black vertex for each vertex and face of $U_{2\eps}$,
and a white vertex for each edge of $U_{2\eps}$ (and $U_\eps$ is the
superposition of $U_{2\eps}$ and its dual). 
However we remove from $U_\eps$ one black vertex $b_0$ on the outer face.
\begin{lemma}
$U_\eps$ has the same number of black and white vertices.
\end{lemma}

This follows simply from the Euler formula $V-E+F=1$ applied to $U_{2\eps}$.
In fact spanning trees of $U_{2\eps}$ rooted at $b_0$ are in bijection
with perfect matchings of $U_\eps$ (see \cite{Temp, KPW}).

Let $G(u,v)$ be the continuous Dirichlet Green's function on $U$.
Let $\tilde G(u,v)$ be the analytic function of $v$ whose
real part is $G(u,v)$. 
On $U_\eps$ we have $K^{-1}=2\eps d\tilde G + o(\eps)$, in the following
sense.
Define $F_0,F_1$ and $F_+,F_-$ by
\begin{eqnarray*}d\tilde G(u,v)&=&\frac12(F_0(u,v)du_x + iF_1(u,v)du_y)\\
&=&\frac14(F_+(u,v)du + F_-(u,v)d\bar u).
\end{eqnarray*}
\begin{thm} When $|w-b|$ is not $O(\eps)$ we have
$$K^{-1}(w,b)=\left\{\begin{array}{llll}
\eps \Re F_0(w,b)+o(\eps)&w\in W_0\mbox{ and }b\in B_0\\
\eps i\Im F_0(w,b)+o(\eps)&w\in W_0\mbox{ and }b\in B_1\\
\eps \Re F_1(w,b)+o(\eps)&w\in W_1\mbox{ and }b\in B_0\\
\eps i\Im F_1(w,b)+o(\eps)&w\in W_1\mbox{ and }b\in B_1.\end{array}\right.$$
\end{thm}

For example on $\R^2$ we have $\tilde G(u,v)=-\frac1{2\pi}\log(v-u)$,
so $d\tilde G= \frac1{2\pi}\frac{du_x+idu_y}{v-u}$ and $F_0=F_1=\frac1\pi
\frac1{v-u}$.
As another example, on the upper half plane we have
$\tilde G(u,v) = -\frac1{2\pi}\log\frac{v-u}{v-\bar u},$ from which we get

$$d\tilde G(u,v) = \frac1{2\pi}\left(\frac{du}{v-u}-\frac{d\bar u}{v-\bar u}
\right).$$
As a consequence 
$F_+(u,v)=\frac2{\pi(v-u)}$ and $F_-(u,v)=-\frac2{\pi(v-\bar u)}.$

It is easier to work with the functions $F_{\pm}=F_0\pm F_1$ since they
transform nicely. 
Recall that $\tilde G$ is conformally invariant, that is if
$\phi:V\to U$ is a conformal homeomorphism then
$\tilde G_V(u,v)=\tilde G_U(\phi(u),\phi(v)).$ 
As a consequence $F_\pm$ are conformally covariant: 
$$F_+^V(u,v)=F_+^U(\phi(u),\phi(v))\phi'(v).$$
$$F_-^V(u,v)=F_-^U(\phi(u),\phi(v))\overline{\phi'(v)}.$$

This allows us to compute asymptotics of $K^{-1}$ on any 
Jordan domain.

More precisely, the result is the following.
When $w,b$ are {\it far apart}, we can use the above
asymptotic formulas. When they are close together, 
the dominant term is controlled by (\ref{close}) with
a deviation of order $\eps$ given by the function(s)
$F_+(u,v) - \frac2{\pi(v-u)}$ and $F_-(u,v)$.

\subsection{Height moment}
Here is a sample computation, the expected value of
$h(p)h(q)$ in the upper half plane. 

\begin{figure}[htbp]
\vskip3in
\PSbox{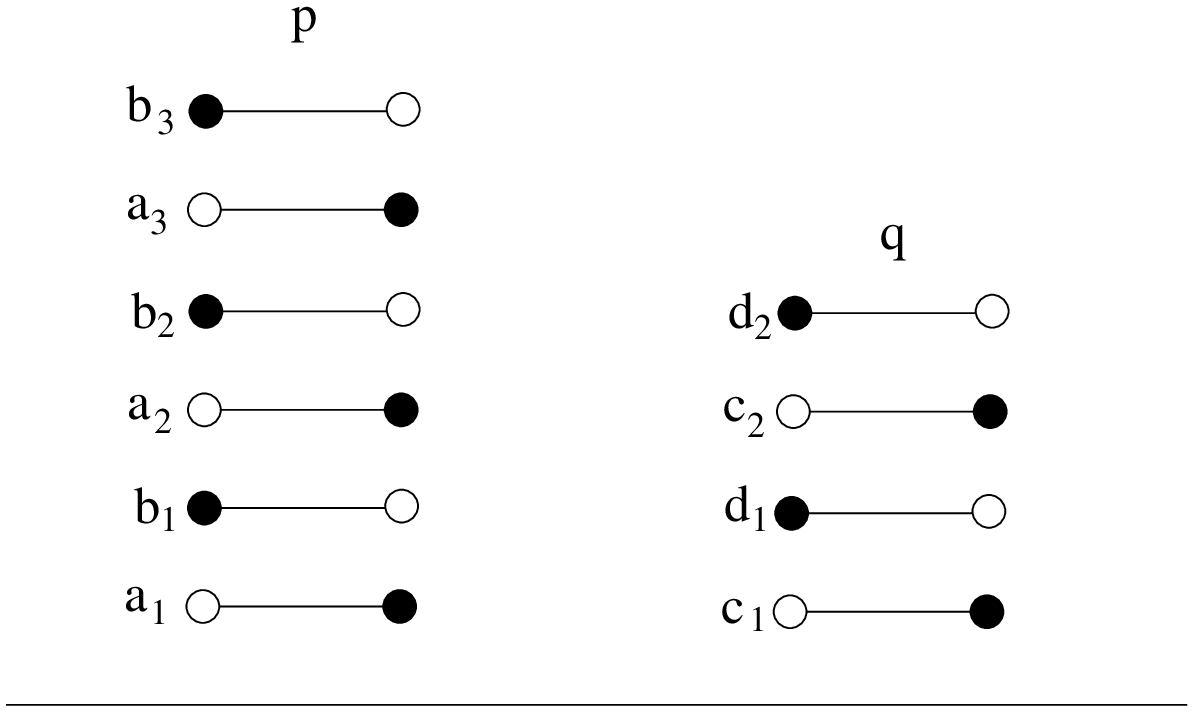}{0in}{0in}
\caption{\label{2mom}}
\end{figure}

Let $a_i,b_i$ be the horizontal edges on a vertical path from $p$ to
the $x$-axis, the $a_i$ being those with a white left vertex,
and similarly let $c_i,d_i$ be those for $q$, see Figure 
\ref{2mom}.
We have $h(p)=4\sum_i a_i-b_i$, $h(q)=4\sum_j c_j-d_j$.
So 
$$\E(h(p)h(q))=16\sum_{i,j}\E(a_ic_j)-\E(b_ic_j)-\E(a_id_j)+\E(b_id_j).$$
Each of these terms can be written as the determinant of
a $2\times 2$ matrix with entries in $K^{-1}$.
For example 
$$\E(a_ic_j)=\det\left(\begin{array}{cc}K^{-1}(a_i^{(0)},a_i^{(1)})&
K^{-1}(c_j^{(0)},a_i^{(1)})\\
K^{-1}(a_i^{(0)},c_j^{(1)}) &K^{-1}(c_j^{(0)},c_j^{(1)})\end{array}\right)$$
where $a_i^{(0)},a_i^{(1)}$ are the white and black vertices of $a_i$,
and likewise for $c_j$.
The diagonal terms will all cancel out, and
after a certain amount of rearranging,  one finds in the limit
$$\E((h(p)-\overline{h(p)})(h(q)-\overline{h(q)}))=\hskip4in$$
$$\hskip1in=\int_{\gamma_1,\gamma_2}
\left|\begin{array}{cc}0&F_+(z_1,z_2)\\F_+(z_2,z_1)&
0\end{array}\right|dz_1dz_2 -\int_{\gamma_1,\gamma_2}
\left|\begin{array}{cc}0&F_-(z_1,z_2)\\\overline{F_-(z_2,z_1)}&
0\end{array}\right|d\overline{z_1}dz_2 -\hskip1in$$
$$\hskip1in\int_{\gamma_1,\gamma_2}
\left|\begin{array}{cc}0&\overline{F_-(z_1,z_2)}\\F_-(z_2,z_1)&
0\end{array}\right|dz_1d\overline{z_2} +
\int_{\gamma_1,\gamma_2}
\left|\begin{array}{cc}0&\overline{F_+(z_1,z_2)}\\\overline{F_+(z_2,z_1)}&
0\end{array}\right|d\overline{z_1}d\overline{z_2}.$$

For the upper half-plane we have $F_+(z_1,z_2)=\frac2{\pi(z_2-z_1)}$ 
and $F_-(z_1,z_2)=-\frac2{\pi(z_2-\overline{z_1})}$.
Plugging these in gives
$$-\frac4{\pi^2}\int_{\gamma_1}\int_{\gamma_2}\frac1{(z_2-z_1)^2}dz_1dz_2 +
\frac4{\pi^2}\int_{\gamma_1}\int_{\gamma_2}\frac1{(z_2-\overline{z_1})^2}
d\overline{z_1}dz_2+
\qquad\qquad$$
$$\qquad\qquad + \frac4{\pi^2}\int_{\gamma_1}\int_{\gamma_2}\frac1{(
\overline{z_2}-z_1)^2} dz_1d\overline{z_2}
-\frac4{\pi^2}\int_{\gamma_1}\int_{\gamma_2}\frac1{(\overline{z_2}-\overline{
z_1})^2}
d\overline{z_1}d\overline{z_2} .$$

The first of these integrals gives
$$-\frac4{\pi^2}\log\frac{(p-q)(r-s)}{(p-s)(r-q)}.$$
Therefore
\begin{eqnarray*}
\E((h(p)-\bar h(p))(h(q)-\bar h(q)))&=&\frac4{\pi^2}\left(
-2\Re\log\frac{(p-q)(r-s)}{(p-s)(r-q)}+2\Re\log\frac{(\overline{p}-q)(r-s)}{
(\overline{p}-s)(r-q)}\right)\\
&=&\frac8{\pi^2}\Re\log\left(\frac{\overline{p}-q}{p-q}\right).
\end{eqnarray*}

Note that this quantity is a multiple of the Dirichlet Green's function
$G(p,q)$.

\subsection{The Gaussian free field}
Let $h_0(p)=h(p)-\overline{h(p)}$ be the fluctuation of $h$ away from the mean.
More generally one finds that all the moments of the height fluctuations
can be written in terms of Green's functions:

\begin{thm}
Let $U$ be a Jordan domain with smooth boundary. Let $p_1,\ldots,p_k\in U$ 
be distinct points. If $k$ is odd we have
$\lim_{\eps\to0}\E(h_0(p_1)\cdots h_0(p_k))=0$. If $k$ is even we have
$$\lim_{\eps\to0}\E(h_0(p_1)\cdots h_0(p_k))=
\left(-\frac{16}{\pi}\right)^{k/2}\sum_{\mbox{pairings } \sigma} 
G(p_{\sigma(1)},p_{\sigma(2)})
\cdots G(p_{\sigma(k-1)},p_{\sigma(k)}).$$
\end{thm}

By Wick's theorem, this implies that $h_0(x)$
is the unique Gaussian process with covariance function
$\E(h_0(x)h_0(y))=G(x,y)$.
This can be taken as a definition of the {\bf massless Gaussian free field}.

An alternate description is as follows.

On $U$ let $\{e_i\}$ be an orthonormal eigenbasis of the Laplacian
with Dirichlet boundary conditions, and $\Delta e_i=\lambda_i e_i$.

Define $$\text{GFF}(x) 
= \sum_{i=1}^\infty \frac{c_i}{\sqrt{|\lambda_i|}}e_i(x)$$
where $c_i$ are independent Gaussians of mean 0 and variance $1$.

This series defines a distribution, not a function: the series diverges
almost surely almost everywhere. However for a smooth function
$\psi$ the series
$$\text{GFF}(\psi)\stackrel{def}{=}
\sum_{i=1}^\infty \frac{c_i}{\sqrt{|\lambda_i|}}\int_U e_i(x)\psi(x)dx$$
converges almost surely.

\begin{thm}[\cite{K.gauss}]As $\eps$ tends to $0$, $h_0$ tends
weakly in distribution to $4/\sqrt{\pi}$ times the massless free field $GFF$ 
on $U$
in the sense that for any smooth function $\psi$ on $U$,
the random variable 
$\eps^2\sum_{x\in U_\eps}\psi(x)h_0(x)$
tends in distribution to $\frac{4}{\sqrt{\pi}}GFF(\psi).$
\end{thm}

It is worth pointing out that, although values of $h_0$ are integral
(after adding a non-random smooth function $\overline{h}$), 
$h_0$ converges to a continuous object.

There is a lot more that can be said about the distribution of
the value of $h$ and $h_0$
at a point, or the joint punctual distributions. 
It can be proved that the distribution of $h$ at a point tends, when scaled
by $\sqrt{\log\frac1{\eps}}$, to a Gaussian with variance
equal to some universal constant.

\section{FK percolation on critical planar graphs}
Up to now we have been working on subgraphs of $\Z^2$. For the dimer model
on certain other regular graphs, such as the honeycomb graph, similar
results can be obtained. However for more general periodic planar graphs,
the situation can be more complicated. In this section we discuss a
different but related model, the FK-percolation model, on an interesting family
of planar graphs the ``isoradial'' graphs. For this family of graphs
we have a surprising property of the partition function and measures
that they depend only on the local structure of the graph, not its
long-range order.

In the subsequent section we return to dimers but again on this same family of
isoradial graphs. 

In this section
we will show that the partition function of the random cluster model has 
a particular form for this special family of planar graphs (which
includes the case $\Z^2$). 

Recall the definition of the FK-percolation (random cluster) model.
Let $G$ be a graph. The space of configurations
is $X=\{w:E\to\{0,1\}\}$ (each edge is open or closed).

Let $\nu:E\to(1,\infty)$ be an assignment of weights to the edges, and
$q>0$ a constant.

We define a probability measure on $X$,
$$\mu(w)=\frac1Z q^c\prod_{e \text{ open}} \nu(e),$$
where $c$ is the number of connected componets of open edges of $w$,
the product is over open edges of $w$, and $Z$ is a normalizing constant,
called the partition function:
$$Z=\sum_{w\in X}q^c \prod_{e \text{ open}} \nu(e).$$

Note that when $q=1$ we are reduced to the standard peroclation model.
More generally when $q$ is an integer the model is equivalent to the 
$q$-state Potts model \cite{Baxter}; in an appropriate
limit $q\to 0$ the model is equivalent to the spanning tree model.

For further information on the model see \cite{Grimm}.

\subsection{Duality} If $G$ is planar let $G^*$ be its planar dual.
Let $\nu^*$ be weights on edges $e^*$ of $G^*$, defined by
$\nu^*(e^*)=q/\nu(e)$. 
\begin{lemma}
For any configuration of edges of $G$ with $c$ components, $k$ edges
and $r$ cycles, we have $c-r+k=V$, where $V$ is the number of vertices of
$G$.
\end{lemma}
Then we have the following
identity between $Z$ and $Z^*$, the dual partition function:
\begin{prop} $$Z^*=Z\cdot q^{E-V}/W,$$
where $W$ is the product of all of the edge weights $W=\prod_{e\in E} \nu(e)$.
\end{prop}

\begin{proof}
Suppose $A$ is a configuration of open edges, with $c$ components
and weight $q^c \prod_{e\in A}\nu(e).$
The dual configuration has edges $E-A$ and $c'$ components.
Its weight is
$$q^{c'}\prod_{e\in E-A}\nu^*(e^*)=q^{c'}\prod_{e\in E-A}\frac{q}{\nu(e)}
=q^{c'}\frac{q^{E-A}}{W}\prod_{e\in A}\nu(e),$$
where $W$ is the product of all edge weights of $G$.
But by the Lemma, $c+A-c'=V$, so the weight of a dual configuration
is a constant $q^{E-V}/W$ times the weight of the primal configuration.
\end{proof}

In particular the duality is an {\bf isomorphism} between the corresponding
probability spaces. As a consequence we have $p_e+p_{e^*}=1$,
where $p_e,p_{e^*}$ are the edge probabilities.

\subsection{$Y-\Delta$ transformation}
\begin{figure}[htbp]
\vskip3in
\PSbox{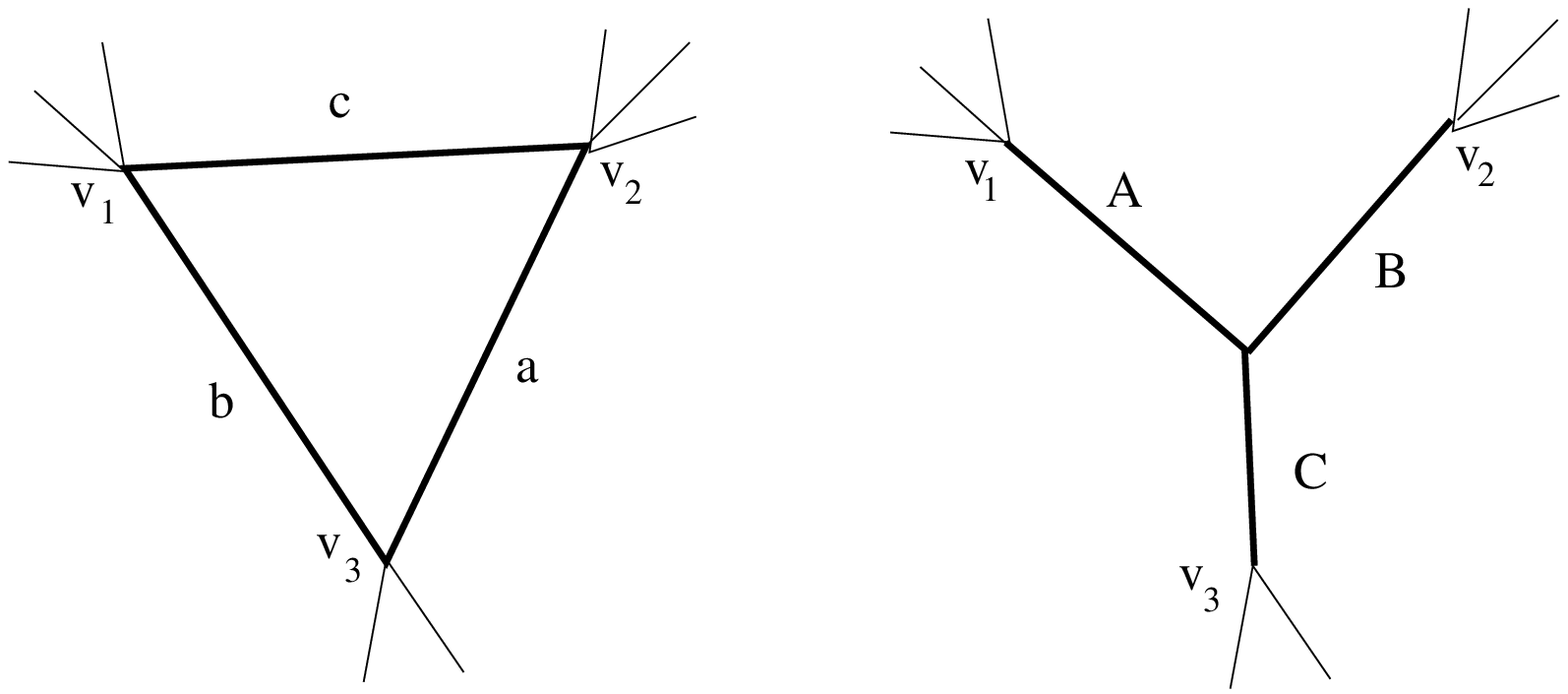}{0in}{0in}
\caption{\label{YD}}
\end{figure}

We can sometimes transform $G$ while preserving the measure $\mu$.
See Figure \ref{YD}.

Given a triangle in $G$
with vertices $v_1,v_2,v_3$ and edge weights $e_{23}=a,
e_{13}=b,e_{12}=c$, and a `$Y$' with the same vertices and edge
weights $A,B,C$, (as in Figure \ref{YD}) we have
\begin{center}
\begin{tabular}{c|c|c}
event&weight in $\Delta$&weight in $Y$\\
\hline
$v_1,v_2,v_3$ connected&$(1+a)(1+b)(1+c)$&$q+A+B+C+AB+BC+AC+ABC$\\
$v_1,v_2$ connected&$(1+c)(1+\frac{a+b+ab}q)$&$q+A+B+C+AB+\frac{AC+BC+ABC}q$\\
$v_1,v_3$ connected&$(1+b)(1+\frac{a+c+ac}q)$&$q+A+B+C+AC+\frac{AB+BC+ABC}q$\\
$v_2,v_3$ connected&$(1+a)(1+\frac{b+c+bc}q)$&$q+A+B+C+BC+\frac{AB+AC+ABC}q$\\
none connected&$1+\frac{a+b+c}q+\frac{ab+ac+bc+abc}{q^2}$&
$q+A+B+C+\frac{AB+AC+BC}q+\frac{ABC}{q^2}$
\end{tabular}\end{center}

Here the events on the left refer to connections via open
edges {\it outside} the triangle.
The central column gives the additional (multiplicative) weight
due to the possible open edges inside the triangle.
The left column gives the additional (multiplicative) weight
due to the possible open edges inside the $Y$.

The $Y-\Delta$ transformation preserves the measures on condition
that the weights in the `$\Delta$' column be proportional
to the weights in the `$Y$' column. This gives $4$ polynomial equations
for the weights $a,b,c,A,B,C$.
\begin{lemma}\label{solnpoly}
These equations have a solution only if
\begin{equation}-q+ab+ac+bc+abc=0\label{critq1}
\end{equation}
and
\begin{equation}
-q^2-q(A+B+C)+ABC=0\label{critq2}.
\end{equation}
In this case the solution is
\begin{equation}A=\frac{q}a,\ \ \ \
B=\frac{q}b\label{aA}\ \ \ \
C=\frac{q}c.
\end{equation}
\end{lemma}

For certain types of graphs/weights, these transformations take on a 
particularly nice form:

\subsection{Isoradial embeddings}
This notion is due to Duffin \cite{Duffin} and Mercat \cite{Mercat}.
We say $G$ has an {\bf isoradial embedding} if $G$ is drawn so that each
face is a cyclic polygon, that is, is inscribable in a circle,
and all the circles have radius $1$.

\begin{figure}[htbp]
\vskip3in
\PSbox{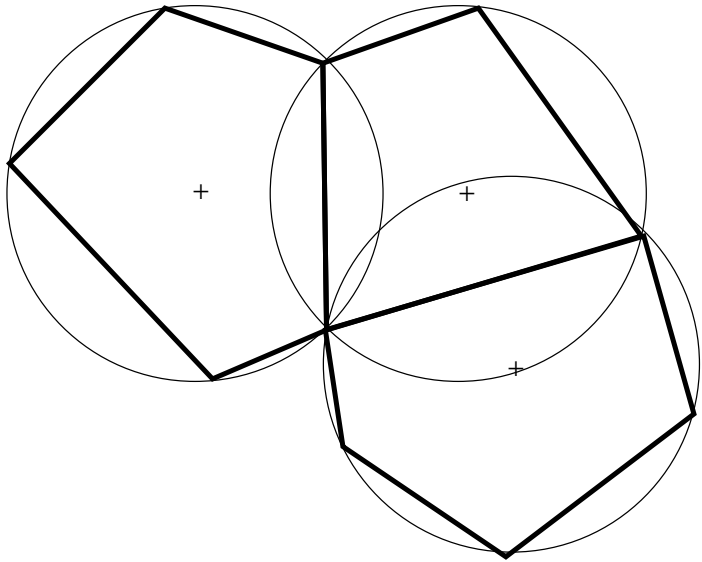}{0in}{0in}
\caption{\label{isorad}Isoradial embedding of a graph}
\end{figure}

\begin{figure}[htbp]
\vskip3in
\PSbox{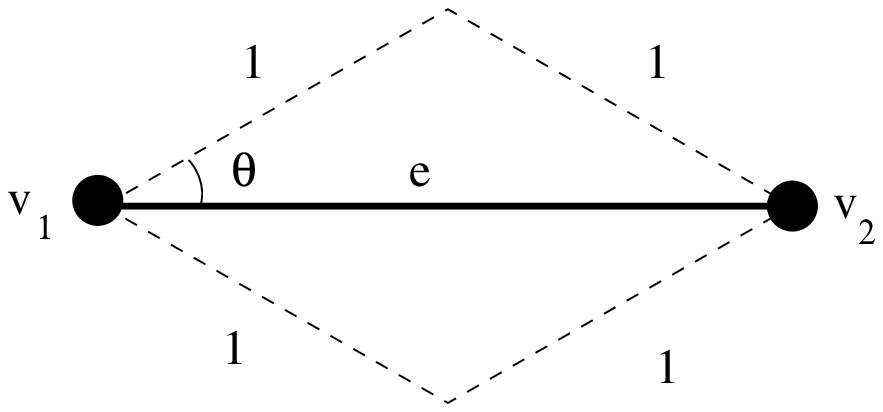}{0in}{0in}
\caption{\label{rhomb}}
\end{figure}

In this case we can embed $G^*$ isoradially as well, using the circle centers
(at least as long as the circle centers are in the interior of the 
corresponding faces of $G$).
Each edge has a rhombus around it whose vertices
are the vertices of the edge and its dual (Figure \ref{rhomb}).

Let $\theta$ be the half-angle of the rhombus and define
$$\nu(e)=\sqrt{q}\frac{\sin(\frac{2r}\pi\theta)}{
\sin(\frac{2r}{\pi}(\frac{\pi}2-\theta))}$$
where $r=\cos^{-1}(\sqrt{q}/2).$

Then one may check that
\begin{lemma}
A $Y-\Delta$ transformation preserves isoradiality under these weights.
\end{lemma}

This follows from plugging in the weights $\nu(e)$ into (\ref{critq1}),
and using the fact that a triangle exists only when
$\theta_a+\theta_b+\theta_c=\pi/2$.

\subsection{Periodic graphs}
A planar
graph $G$ is {\bf periodic} if translations in $\Z^2$ are (weight-preserving)
isomorphisms.  In this case $G_n=G/n\Z^2$ is a toroidal graph.
Let $Z_n$ be the partition function on $G_n$.
we define $Z = \lim_{n\to\infty} Z_n^{1/|G_n|}.$
This is the {\bf partition function per site} for the infinite
periodic graph $G$.

\subsection{Main theorem}
Here is the main result of this section:
\begin{thm}
There exists a function $F_q:[0,\frac{\pi}2]\to\R$ such that
for any periodic, isoradial planar graph $G$, with edge weights
$\nu$ as above, the partition function per site is
$$\log Z = -\frac12\log q + \frac1{|G_1|}
\sum_{\text{edges in a f.d.}} F_q(\theta).$$
\end{thm}

The proof doesn't tell us anything about the function $F_q$.
However $F_q$ was computed (non-rigorously) by Baxter \cite{Baxter}
(for the graph $\Z^2$, and hence for any isoradial graph).
For $q\in(0,4)$ the answer is
$$F_q(\theta)= -\frac12\int_{-\infty}^{\infty}\frac{\sinh((\pi-r)t)\sinh(
\frac{4r}{\pi} \theta t)}{t\sinh(\pi t)\cosh(r t)}dt.$$
Any information about this integral would be greatly appreciated.
\medskip

\noindent{\bf Proof sketch:}
Each $Y\leftrightarrow\Delta$ changes the partition function by $$Z_Y=Z_\Delta 
\frac{q^2}{abc}.$$ Define 
a new-fangled ``partition function''
$$\tilde Z=\frac{Zq^{V/2}}{\sqrt{\prod_E\nu(e)}}.$$
It is easy to check that $\tilde Z_Y=\tilde Z_\Delta$. 

\begin{figure}[htbp]
\vskip3in
\PSbox{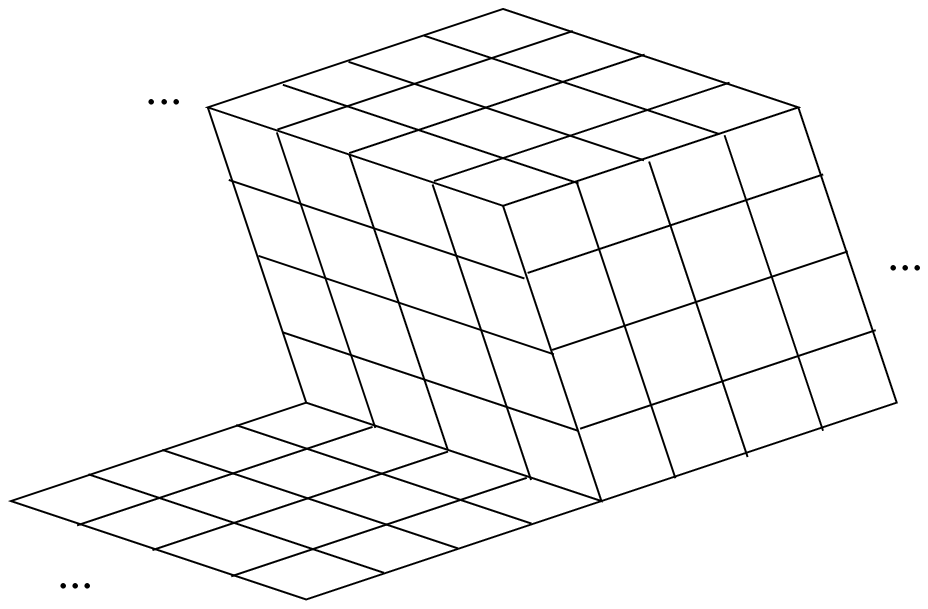}{0in}{0in}
\caption{\label{blocks}}
\end{figure}

The idea is to use $Y-\Delta$ transformations to turn the graph
into a graph with large blocks which are copies of big pieces of $\Z^2$.
Because $G$ is periodic, each bi-infinite rhombus chain is parallel 
to copies of itself. We'll rearrange the graph so as to move these
chains adjacent to each other in clumps, so that the resulting
rhombus tiling consists of large rhombi each tiled with $N^2$ copies
of smaller versions of themselves (Figure \ref{blocks}).

\begin{figure}[htbp]
\vskip3in
\PSbox{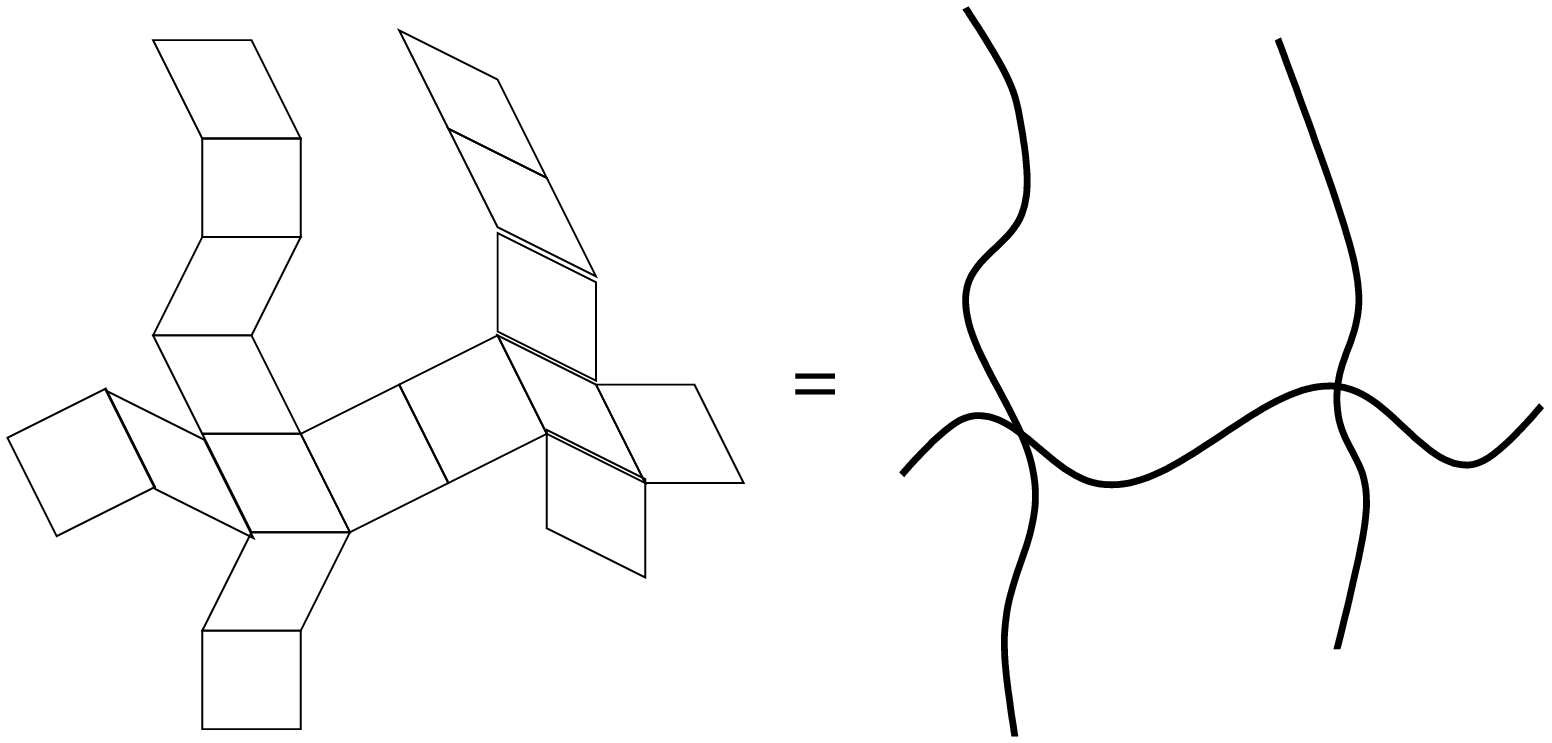}{0in}{0in}
\caption{\label{strings}}
\end{figure}

How do we do this transformation?
Convert rhombus chains into strings as in Figure \ref{strings}

Slide the strings around without creating new intersections. When a string
crosses the intersection of another pair of strings, the
crossing corresponds to a $Y-\Delta$ transformation on the underlying graph.
We can slide the strings around 
until the graph looks like a union of blocks, each tiled
by parallel rhombi, as in Figure \ref{blocks}.
That is, we can do this unless some rhombus chain is separated from its
translates by a third chain which does not intersect either but
whose rhombi have a different common parallel. In this case
it is possible to show that simply exchanging the chains does not change the 
partition function. (One way to see this is to add an extra edge and it's
``negative'' at the same place so that you can exchange the two rows using
a sequence of $Y-\Delta$ moves...)

Now by subadditivity, when there are only large blocks, the partition
function can be computed separately on each block and the results multiplied
together (with an error of lower order). But on each block the graph is
a piece of $\Z^2$ and the partition function is $\exp(Ac(1+o(1)))$
where $A$ is the area and $c$ a function of the angle only.
Letting $N\to\infty$ we can conclude in that
$\tilde Z=\prod_{\text edges}F_q(\theta)$ for some $F_q$.
The result follows.  \hspace*{\fill} $\square$\medskip 

This theorem has some nontrivial consquences even without knowing $F_q$.
For example let $Z_{rect}$ be the partition function per site
for the graph $\Z^2$ with edge lengths $1$ vertically and $\sqrt{3}$
horizontally (and corresponding weights $\nu$ as above). 
Let $Z_{hex}$ be the partition function per site
for the honeycomb graph with edge lengths $1$.
Then
$$\frac{Z_{\text{rect}}}{Z_{\text{hex}}}=\frac{2\cos(\frac23\cos^{-1}(
\frac{ \sqrt{q}}{2}))}{q^{1/6}}.$$

As another ``application'' of the theory, it is natural to conjecture that
the FK-percolation model is critical for weights $\nu$.
For example in standard percolation ($q=1$) we should have
$p_c/(1-p_c)=\nu$. 

For example on the
on the honeycomb lattice the critical probability should be
$$p_c= \frac{2\cos\frac{\pi}9}{1+2\cos\frac{\pi}9}\approx 0.652704.$$
This was in fact proved for the honeycomb lattice by Wierman \cite{W}.

\begin{prop}The only function $f$ for which $(a,b,c)=(f(\theta_a),
f(\theta_b),f(\theta_c))$ satisfies (\ref{critq1}) when
$\theta_a+\theta_b+\theta_c=\pi/2$ and (\ref{aA}) is
$$f(\theta) = \sqrt{q}\frac{\sin\left(\frac{2r}\pi\theta \right)}{
\sin\left(\frac{2r}\pi(\frac{\pi}2-\theta)\right)}$$
where $r=\cos^{-1}(\sqrt{q}/2).$
\end{prop}

\section{Integrability and dimers on critical planar graphs}
Surprisingly, the same isoradiality condition which worked so well
in the random cluster model also has a simplifying effect on the
dimer partition function.

Let $G$ be a bipartite planar graph.
We embed $G$ isoradially. Let $\nu(e)=\sin2\theta$,
where $\theta$ is the half-angle of the corresponding rhombus.
We use these weights to define a probability measure on perfect matchings
of (finite subgraphs of) $G$ as we did before: the probability of
a matching is proportional to the product of its edge weights.

We define a Kasteleyn matrix as follows:
if $e=wb$ is an edge with dual edge $e^*=p^*q^*$,
\begin{figure}[htbp]
\vskip3in
\PSbox{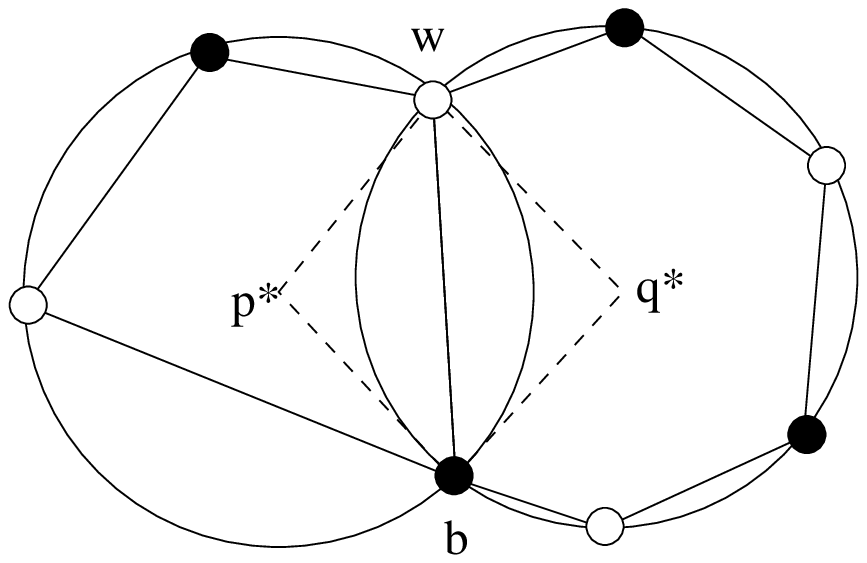}{0in}{0in}
\caption{\label{isoraddimer}}
\end{figure}
then $K(w,b)= i(p^*-q^*).$
\begin{lemma}
$K$ is a Kasteleyn matrix for (any finite simply connected subpiece of) $G$.
\end{lemma}

Note that $|K(w,b)|=2\sin\theta$.

\begin{proof}
It suffices to show that $K$ is Kasteleyn-flat
\cite{Kuper},
that is, we must show that for each face of $G_1$ with
vertices $u_1,v_1,\dots,u_m,v_m$ in cyclic order, we have
$$\arg(K(u_1,v_1)\dots K(u_m,v_m))=\arg(
(-1)^{m-1}K(v_1,u_2)\dots K(v_{m-1},u_m)K(v_m,u_1)).$$
(This identity implies that two dimer configurations which only differ
around a single face have the same argument in the expansion of
the determinant.
By \cite{Kuper}, any two configurations can be obtained from one another
by such displacements.)
To prove this identity, note that it is true if the points are
regularly spaced along a $2m$-gon, and note that it remains
true if you move one point at a time.
\end{proof}

There are two main results in this section, an explicit
computation of $Z$ and $K^{-1}$.
\begin{thm}
The determinant per site of $K$ satisfies
$$\log Z=
\frac1N\sum_{\text{edges } e}\frac1{\pi}L(\theta) + \frac{\theta}{\pi}\log
2\sin\theta,$$
where $\nu(e)=2\sin(\theta(e))$,
$N$ is the number of vertices in a fundamental domain, the sum is over the
edges in a fundamental domain, and $L$ is the Lobachevsky function,
$$\label{Lob}L(x)=-\int_0^x\log2\sin tdt.$$
\end{thm}

To describe $K^{-1}$, first define for any complex parameter $z$
an {\bf elementary discrete analytic function} to be a function
satisfying
$f_{v'}=f_v/(z-e^{i\alpha})$, if $vv'$ leads away from a white, or towards
a black vertex, and 
$f_{v'}=f_v\cdot (z-e^{i\alpha})$, if $vv'$ leads away from a black, or towards
a white vertex.
Such an $f$ is well-defined the product of the multipliers
going around a rhombus leads back to the starting value.
\begin{thm}
We have
\begin{equation}
\label{Kinverse}
K^{-1}(b,w_0)=
-\frac{1}{2\pi i}\sum_{\text{poles } e^{i\theta}} \theta\cdot
\mbox{Res}_{z=e^{i\theta}}(f_b/f_{w_0}),
\end{equation}
where the angles $\theta\in\R$ are chosen appropriately.
This can be written
$$\frac1{4\pi^2i}\int_C \frac{f_b(z)}{f_{w_0}(z)}\log z dz,$$
where $C$ is a closed contour surrounding cclw the part of the circle
$\{e^{i\theta}~|~\theta\in[\theta_0-\pi+\eps,\theta_0+\pi-\eps]\}$
which contains all the poles of $f_b$,
and with the origin in its exterior.
\end{thm}

\begin{proof}
Let $F(b)$ denote the right-hand side of (\ref{Kinverse}).
We will show that $\sum_{b\in B}K(w,b)F(b) = \delta_{w_0}(w)$, and
$F(b)$ tends to zero when $b\to\infty$.

Let $C_{\theta}$ be a small loop around $e^{i\theta}$ in $\C$.
We have
\begin{eqnarray*}
\sum_{b\in B}K(w,b)F(b)
&=& \sum_{j=1}^ki(-e^{i\theta_j}+e^{i\theta_{j-1}})F(b_j),\\
&=&-\frac{1}{2\pi i}\sum_{j=1}^ki(-e^{i\theta_j}+e^{i\theta_{j-1}})
\sum_{\text{poles }e^{i\theta}}
\theta\cdot \text{Res}_{z=e^{i\theta}}(f_{b_j})\\
&=&-\frac1{2\pi}\sum_{\text{poles }e^{i\theta}}\theta\cdot\sum_{j=1}^k
(-e^{i\theta_j}+e^{i\theta_{j-1}})\frac1{2\pi i}\int_{C_\theta} f_{b_j}(z)dz\\
&=&-\frac1{2\pi}\sum_{\text{poles }e^{i\theta}}\frac{\theta}{2\pi i}
\int_{C_\theta} f_{w}(z)
\sum_{j=1}^k
\frac{(-e^{i\theta_j}+e^{i\theta_{j-1}})}{(z-e^{i\theta_j})(z-e^{i\theta_{j-1}})}dz\\
&=&-\frac1{2\pi}\sum_{\text{poles }e^{i\theta}}\frac{\theta}{2\pi i}
\int_{C_\theta} f_{w}(z)\left(\sum_{j=1}^k
\frac1{(z-e^{i\theta_{j-1}})}-\frac1{(z-e^{i\theta_j})}\right)dz\\
&=&-
\frac1{2\pi}\sum_{\text{poles }e^{i\theta}}\frac{\theta}{2\pi i}
\int_{C_\theta} 0 dz = 0.
\end{eqnarray*}

However when $w=w_0$ we have
$$F(b_j)=-\frac1{2\pi i}\left(\frac{\theta_j-\theta_{j-1}}{e^{i\theta_j}-e^{i\theta_{j-1}}}
\right),$$
so that
$$\sum_{j=1}^ki(-e^{i\theta_j}+e^{i\theta_{j-1}})F(b_j)=\frac1{2\pi}\sum_{j=1}^k
\theta_j-\theta_{j-1} = 1,$$
since the angles increase by $2\pi$ around $w_0$.

To complete the proof, one must show that $F(b)\to0$ as $b\to\infty$.
\end{proof}

\begin{figure}[htbp]
\vskip3in
\PSbox{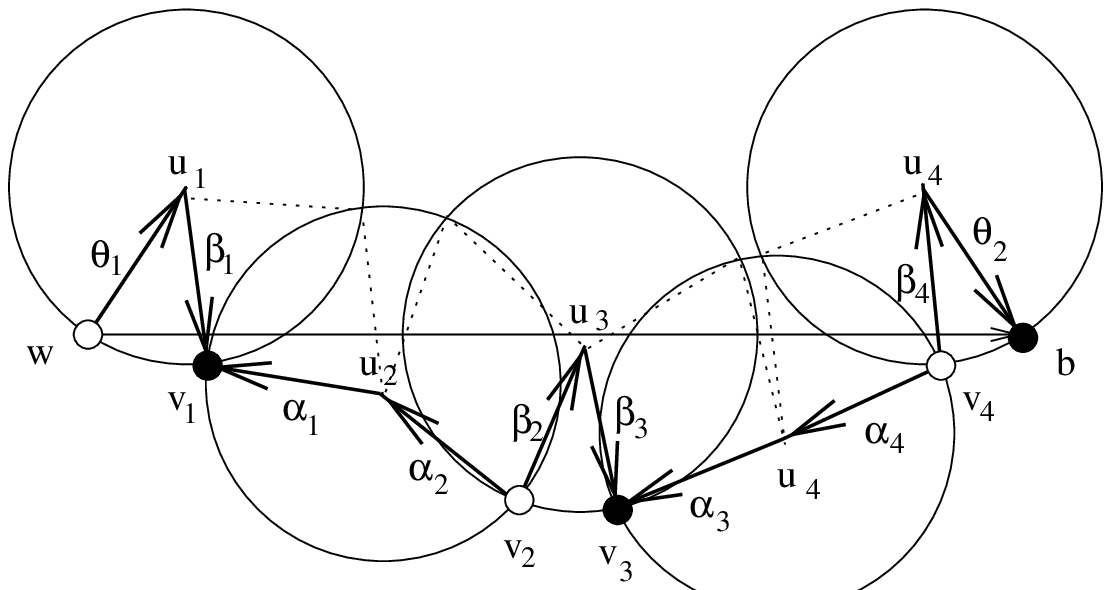}{0in}{0in}
\caption{\label{circlepath}}
\end{figure}

\begin{prop}\label{1/z}
We have
$$K^{-1}(b,w) = \frac{1}{2\pi}\left(\frac{1}{b-w}+
\frac{\gamma}{(\bar b-\bar w)} \right)+
O\left(\frac1{|b-w|^2}\right),$$
where $\gamma=f_b(0)=e^{-i(\theta_1+\theta_2)}\prod e^{i(-\beta_j+\alpha_j)},$
and the angles are as illustrated in Figure \ref{circlepath}.
\end{prop}

\subsection{determinants}
For a finite graph $G$ we define the {\bf normalized determinant},
or determinant {\bf per site}, 
$\det_1 M$
of an operator $M\colon\C^{G}\to\C^{G}$ to be
$$\mbox{det}_1 M\stackrel{\text{def}}{=}|\det M|^{1/|G|},$$ where
$|G|$ is the number of vertices of $G$.

For an operator $M$ on a finite graph,
$\det M$ is a function of the matrix entries $M(i,j)$
which is linear in each entry separately.
In particular for an edge $e=ij$,
as a function of the matrix entry
$M(i,j)$ we have $\det M=\alpha+\beta M(i,j)$,
where $\beta$ is $(-1)^{i+j}$ times
the determinant of the minor obtained by removing
row $i$ and column $j$.
That is,
$$\frac{\partial(\det M)}{\partial(M(i,j))}= M^{-1}(j,i)\cdot\det M,$$
or
\begin{equation}\label{detdef}
\frac{\partial(\log\det M)}{\partial M(i,j)}= M^{-1}(j,i).
\end{equation}

Suppose that $G$ is periodic
under translates by a lattice $\Lambda$.
Let $G_n = G/n\Lambda$, the finite graph which is the
quotient of $G$ by $n\Lambda$. Now if we sum (\ref{detdef}) for all
$\Lambda$-translates of edge $ij$ in $G_n$, and then divide
both sides by $|G_n|$, it yields
\begin{equation}\label{det1def}
\frac{\partial(\log\det_1 M)}{\partial M(i,j)}= \frac1{|G_1|}M^{-1}(j,i),
\end{equation}
where $M(i,j)$ is now the common weight of all translates of edge $ij$,
that is, the left-hand side is the change in $\log\det_1 M$ when
the weight of all translates of $ij$ changes.
Note that this equation is independent of $n$.

It is now a short computation to compute the determinant of $K$
from the exact form of $K^{-1}$....see \cite{K.isorad} 

\subsection{Isoradial embeddings}
The set of isoradial embeddings $\text{ISO}(G)$, 
when parametrised by the rhombus angles,
is convex. This follows because an isoradial embedding is determined
from a set of rhombus angles by the linear conditions that the sum
around each vertex must be $\pi$.

There is a unique (possibly degenerate)
isoradial embedding maximizing $Z$, because $Z$ is
strictly concave on $\text{ISO}(G).$

\begin{figure}[htbp]
\vskip3in
\PSbox{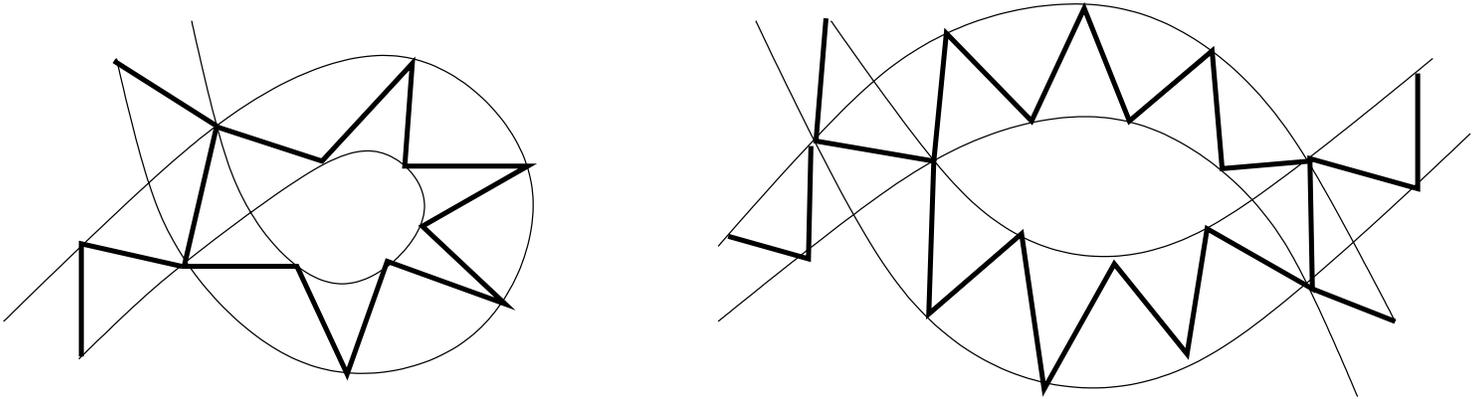}{0in}{0in}
\caption{Not allowed in isoradial graphs\label{zigzag}}
\end{figure}

(Joint with J-M. Schlenker) A {\bf zig-zag} path in a planar graph
is a path which turns maximally left at a vertex, 
then maximally right at the next vertex, then maximally left, and so on.
That is, it leaves from a vertex along an edge which is adjacent in cyclic
order to the edge it entered on, alteranting to the right and to the left.
A planar graph has an isoradial embedding
if and only if the following two conditions are satisfied 
(Figure \ref{zigzag}): 

1. No zig-zag path crosses itself 

2. No two zig-zag paths cross more than once.

It is not hard to see why these conditions are necessary: 
a zig-zag path corresponds to a rhombus chain in an isoradial 
embedding, and such a chain is monotone and so cannot cross itself.
Two such chains cannot cross more than once due to the orientation
of their common rhombi.

\subsection{Hyperbolic ideal polyhedra}
To an isoradially embedded graph $G$ we can associate an ideal
polyhedron in $\H^3$ as follows. Suppose $G$ is embedded
on the $xy$ plane which is the boundary at $\infty$ of the upper
half-space model of $\H^3$. Vertices of $P$ are vertices of 
$G^*$ (circle centers) and edges of $P$ are geodesics which project
vertically to the edges of $G^*$. 

The volume of $P$ per fundamental domain can be seen to be the same as the
entropy per fundamental domain of the dimer model on $G$.
That is, 
$$\frac1\pi\mbox{Vol}(P)=\frac1\pi\sum_{e\in \text{f.d.}} 
L(\theta) = \log Z - \sum_{e\in \text{f.d.}} \Pr(e)\log\nu(e).$$

Moreover the term 
$$\sum_{e\in \text{f.d.}} \Pr(e)\log\nu(e)=\sum_{e\in \text{f.d.}}
\frac{\theta}{\pi}\log2\sin\theta$$
can be associated with the mean curvature of $P$,
which is by definition
the sum over the edges of the dihedral angle times the (normalized)
hyperbolic edge length. 
\medskip

\noindent{\bf Question:} 
Can one construct a tiling in some canonical (entropy-preserving)
way using the geodesic flow inside $P$?

\end{document}